\documentclass[reqno, 12pt, draft]{amsart} % 
\usepackage{amssymb, graphicx} 
\theoremstyle{plain} % 
\newtheorem{theorem}{\indent\sc Theorem}[section] % 
\newtheorem{lemma}[theorem]{\indent\sc Lemma}
\newtheorem{corollary}[theorem]{\indent\sc Corollary}
\newtheorem{proposition}[theorem]{\indent\sc Proposition}

\theoremstyle{definition} % 
\newtheorem{definition}[theorem]{\indent\sc Definition}
\newtheorem{remark}[theorem]{\indent\sc Remark}

\newtheorem{problem}[theorem]{\indent\sc Problem}

\begin{document}

\title[Equivalence problem and double fibration]{\bf{Equivalence Problem for Second Order PDE and Double Fibration as a Flat Model Space}} 

\author[T. Noda]{Takahiro Noda} % 

\dedicatory{Dedicated to Professor Hajime Sato on his first retirement}


\address{ 
Takahiro Noda \endgraf
Graduate School of Mathematics \endgraf
Nagoya University \endgraf
Chikusa-ku, Nagoya 464-8602 \endgraf
Japan
}
\email{m04031x@math.nagoya-u.ac.jp}

%%%%%%%%%%%%%%%%%%%%%%%%%%%%%%%%%%%%%%%%%%%%%%%%%%%%%%%

\maketitle
\begin{abstract}
In this paper, we consider an equivalence problem of second order partially differential
equations (PDE) and a duality of the flat differential equation. 
For the equivalence problem, explicit form of invariants (curvatures) are given.
 In particular, if all of the curvatures vanish, 
then PDE are equivalent to the flat equation. 
We also investigate a duality associated with the flat equation
using double fibrations. These double fibrations are described in terms of transformation
 groups.
\vspace{0.5cm}

\end{abstract}

\section{Introduction}

\hspace{4mm} 
In this paper, we investigate second order PDE for one unknown function 
of two variables. That is, we consider a problem for when these equations 
are equivalent to the flat equation, and we also consider a duality for the flat equation. 
The equivalence problem for differential equations is simply explained as follows.
We fix classes of differential equations and a group of coordinate transformations. 
Then, we consider a problem how differential equations change under 
coordinate transformations. 
We can also express this problem in terms of group actions. Let $\mathcal G$ be a coordinate transformation group and $X$ be a set of certain differential equations. 
Then the equivalence problem for differential equations in $X$ is interpreted as the 
problem of determining the orbit decomposition with respect to the action of 
$\mathcal G$ on $X$. 
 
The equivalence problem is studied deeply by Sophus Lie and $\acute{\rm E}$lie Cartan, 
and many other authors. We mention a few historical background here. 
(See [10] for a detailed history of the equivalence problem.)
Sophus Lie studied an action of the contact diffeomorphism group 
$\mathcal G:={\rm Cont(\mathbb R^3)}$ on $X:=\left\{y''=f(x,y,y')\right\}$, 
and obtained the fact that this action is transitive. In this case, the
orbit decomposition of $X$ for the action of $\mathcal G$ has just one orbit. 
After the work of S.Lie, A.Tresse studied the following case. 
Let $\mathcal G$ be the subgroup ${\rm Diff(\mathbb R^2)^{cont}}$ consisting of lifts
of diffeomorphisms on $\mathbb R^2$ to the jet space $J^{1}(\mathbb R, \mathbb R)$,
and same set of differential equations $X:=\left\{y''=f(x,y,y')\right\}$ .
Under this setting up, Tresse considered a orbit decomposition of the 
action of $\mathcal G$ on $X$. Contrary to the above problem considered by Lie, Tresse proved that this action is not transitive.\par
 On the other hand, $\acute{\rm E}$lie Cartan
 also considered the same problem from a different method, which is now called the equivalence method ([3], [12], [18]).\par
 Along this historical background, we consider an equivalence problem for 
second order PDE for one unknown function of two variables $y=y(x_1,x_2)$:
\begin{equation}
\frac{\partial ^{2} y}{\partial x_i \partial x_j}=f_{ij}(x_{1},x_{2},y,z_{1},z_{2}),
\end{equation}
where, $f_{ij}$ ($1 \leq i,\ j\leq 2$) satisfying 
$f_{ij}=f_{ji}$ are $C^{\infty}$ functions on 
$J^{1}(\mathbb R^{2}, \mathbb R):={\bigl\{(x_{1}, x_{2}, y, z_{1}, z_{2})}\bigr\}$, and 
$z_{1}=y_{x_{1}}$, $z_{2}=y_{x_{2}}$. 
If $f_{ij}$ all vanish, (1) is called the flat equation. 
We take the group ${\rm ScaleDiff(\mathbb R^3)}^{\rm cont}$ of lifts of scale transformations 
on $\mathbb R^3$ as a transformation group $\mathcal G$.\par
 We will calculate explicitly the curvatures for this equivalence problem by using 
Cartan's equivalence method. We obtain the necessary and sufficient condition when the 
second order PDE satisfying integrability condition is equivalent to 
the flat equation via a vanishing condition of these curvatures ([18]). 
Then, our main theorem can be stated as follows.
\vspace{0.5cm}

{\bf \sc Main Theorem 1.}
{\it
For the above equivalence problem, we determine the fifteen curvatures $M_i,\ S_j$ 
explicitly $($curvatures $M_i,\ S_j$ are given in page $11)$. 
In particular, we consider the equation {\rm(1)} for the following functions $f_{ij}:$
\begin{center}
$f_{11}=P(x_1,x_2,y),\hspace{1cm} f_{12}=Q(x_1,x_2,y),\hspace{1cm} f_{22}=R(x_1,x_2,y).$
\end{center}
Then, this equation is {\rm(}locally{\rm)} equivalent to the flat equation 
under lifts of scale transformations if and only if this equation is integrable.
}
\vspace{0.5cm} 

Compare with equivalence problems of second order ODEs, there is a lot of
curvatures in this theorem. 
The reason is given by the following consideration. 
In general, orbit decompositions for PDEs are more complicated than orbit decompositions 
for ODEs. Moreover, $\mathcal G={\rm ScaleDiff(\mathbb R^3)}^{\rm cont}$ is a very strongly 
restricted group. Therefore, this result is obtained. 
 Conversely, if we take groups larger 
than ${\rm ScaleDiff(\mathbb R^3)}^{\rm cont}$, then we obtain a few of curvatures. 
For example, we can consider the group $\mathcal G={\rm Diff(\mathbb R^3)}^{\rm cont}$
as a such group.\par
We also discuss a duality associated with differential equations 
via double fibration. In particular, we consider a 
duality between the coordinate space and the solution space of the flat equation. 
Double fibrations play an important role for a study of this duality. Moreover, 
these fibrations are usually described via some transformation groups appeared in 
equivalence problem ([2], [12]). 
For the group ${\rm ScaleDiff(\mathbb R^3)^{cont}}$, we can not obtain a fibration of 
compact type, because the group ${\rm ScaleDiff(\mathbb R^3)^{cont}}$ is 
too small. Hence, it is natural to consider an existence problem of groups from which 
double fibration of compact-type is obtained as a flat model space. For this problem, 
we find a non-trivial group which gives a fibration of compact-type: 
\begin{align}\label{compact}
G&=\left\{g\in SL(4,\mathbb R) \ |\ g[e_3]=[e_3],\ ^{t}g^{-1}[e_3]=[e_3] \right\} 
\nonumber \\
 &=\left\{
\begin{pmatrix}
*&*&0&* \\
*&*&0&* \\
0&0&*&0 \\
*&*&0&* 
\end{pmatrix}
\in SL(4, \mathbb R) \right\} \nonumber
\end{align}
For this group, we obtain the following fibration of compact-type.

{\sc Main Theorem 2.}
{\it
A double fibration constructed by the above group is the following 
fibration of compact-type.
\begin{equation*}
\setlength{\unitlength}{0.7mm}
\begin{picture}(100,50)(-30,0)
\put(3,0){\makebox(-20,10)[r]{$G/(G\cap H_4)\cong  RP^2$}}
\put(16,30){\makebox(-10,10){$G/(G\cap H)\cong \mathbb F_{V^3}(1,2)$}}
\put(39,0){\makebox(20,10)[l]{$G/(G\cap \overline H_1)\cong RP^2$}}
\put(-15,17){\makebox(-50,10){$(G\cap H_4)/(G\cap H)\cong S^1$}}
\put(50,16){\makebox(35,10){$(G\cap \overline H_1)/(G\cap H)\cong S^1$}}
\put(1,30){\vector(-1,-1){20}}
\put(20,30){\vector(1,-1){20}}
\end{picture}
\end{equation*}
}
The coordinate transformation group $\mathcal G$ corresponding to this group $G$ is constructed by the transformations of the form:
$$X_1=X_1(x_1,x_2),\ 
X_2=X_2(x_1,x_2),\ 
Y=\frac{y}{A(x_1,x_2)}.$$

Finally, we will consider the dual equations of original equations and calculate 
these dual equations explicitly ([17]).\par 
\hspace{-0.5cm}$\bf{Acknowledgements.}$
The author would like to thank Professors. Hajime Sato and\\ Tatsuya Tate for their 
lectures and supports through this work. 
\section{Equivalence problem and $G$-structure}
 \hspace{4mm}
 In this section, we introduce an equivalence problem and explain the $G$-structure 
associated with this problem. For this purpose, we prepare some terminology and
 notation.\par
For functions of two variables $y=y(x_1,x_2)$, we consider the second order PDE (1), 
and diffeomorphisms $\phi$ on $\mathbb R^3$ of the form 
\begin{equation*} \phi(x_{1}, x_{2}, y)=(X_{1}(x_{1}), X_{2}(x_{2}), Y(x_{1}, x_{2}, y)).
\end{equation*} 
The map $\phi$ of this form is called a scale transformation. 
A scale transformation $\phi$ 
lifts naturally to a contact diffeomorphism $\hat{\phi}$ of 
$J^1(\mathbb R^2, \mathbb R)$ defined by:
\begin{equation*} 
{ \hat\phi(x_1, x_2, y, z_1, z_2)=(X_{1}(x_1), X_{2}({x_2}), Y(x_{1}, x_{2}, y), Z_1, Z_2)},
\end{equation*}
where, $Z_1=\frac{Y_{x_{1}}+Y_{y}z_{1}}{(X_{1})_{x_1}}$, 
$Z_2=\frac{Y_{x_{2}}+Y_{y}z_{2}}{(X_{2})_{x_2}}$.
We can easily check that the map $\hat{\phi}$ is a contact diffeomorphism:
\begin{equation*} 
\begin{split} 
\hat\phi^{*}(dy-z_{1}dx_{1}-z_{2}dx_{2})   &=dY-Z_{1}dX_{1}-Z_{2}dX_{2} \\
                                          &=Y_{x_1}dx_1+Y_{x_2}dx_2+Y_{y}dy \\
  &-\frac{Y_{x_{1}}+Y_{y}z_{1}}{(X_{1})_{x_1}}(X_{1})_{x_1}dx_1 
  -\frac{Y_{x_{2}}+Y_{y}z_{2}}{(X_{2})_{x_2}}(X_{2})_{x_2}dx_2 \\
   &=Y_{y}(dy-z_{1}dx_{1}-z_{2}dx_{2}).
\end{split} 
\end{equation*}
We introduce the following terminology:
\begin{align*}
{\rm ScaleDiff(\mathbb R^3)}:&=\left\{\rm Scale\ transformation\ on\ \mathbb R^3 \right\}, \\
{\rm Diff(\mathbb R^3)^{cont}}:
&=\left\{{\rm The\ lift\ of\ \rm Diff(\mathbb R^3)\ to\ J^{1}(\mathbb R^{2}, \mathbb R)}\right\}, \\{\rm ScaleDiff(\mathbb R^3)^{cont}}:
&=\left\{{\rm The\ lift\ of\ \rm ScaleDiff(\mathbb R^3)\ 
to\ J^{1}(\mathbb R^{2},\mathbb R)}\right\}, \\
X:&=\left\{\rm second\ order\ PDE\ (1)\right\}.
\end{align*}
The main problem in the present paper is the following.
\begin{problem}
Examine the orbit decomposition under the action of\\ 
$\rm ScaleDiff(\mathbb R^3)^{cont}$ on $X$. 
\end{problem}
In order to resolve the above problem, 
we use a $G$-structure associated with the equation (1).
First, we replace from data of second order PDE (1) to data of differential system
([3], [12], [19]).
We choose the following coframe of $J^{1}(\mathbb R^{2},\mathbb R)$ corresponding to the 
equation (1),
\begin{align}
\underline{\theta}_{0}:&=dy-z_{1}dx_{1}-z_{2}dx_{2}, \nonumber \\
\underline{\theta}_{1}:&=dz_1-f_{11}dx_{1}-f_{12}dx_{2}, \nonumber \\
\underline{\theta}_{2}:&=dz_2-f_{21}dx_{1}-f_{22}dx_{2}, \\
\underline{\omega}_{1}:&=dx_1, \nonumber \\
\underline{\omega}_{2}:&=dx_2. \nonumber
\end{align}
We consider the Frobenius system 
\begin{equation}
\mathcal{I}:=\bigl\{\underline{\theta}_{0}, 
\underline{\theta}_{1}, \underline{\theta}_{2}\bigr\}_\text{ diff } \text{ with }
 \underline{\omega}_{1}\wedge\underline{\omega}_{2} \not =0
\end{equation}
 constructed by this coframe. 
The correspondence between second order PDE (1) and the Frobenius system $\mathcal{I}$ 
is described as follows. 
Consider vector fields on $J^{1}(\mathbb R^{2},\mathbb R)$ which 
are annialated by $\underline{\theta}_{i}$, while are not annialated by 
$\underline{\omega}_{i}$. At any point on $J^{1}(\mathbb R^{2},\mathbb R)$, 
such vector fields are generated by two vector fields $v_1,\ v_2$. The integral surfaces 
which are tangent to the 2-plane $span\bigl\{v_1,\ v_2\bigr\}$ at any point are the graphs 
of solutions of the second order PDE (1). 
Then, the parameters $(x_1$, $x_2)$ are regarded as a local coordinate system of this integral surface.\par
The Frobenius condition (integrability condition) of the Frobenius
system $\mathcal{I}$ is:
\begin{equation*}
d\underline{\theta}_i \equiv 0 \hspace{5mm} 
(\text{mod } \underline{\theta}_{0}, \underline{\theta}_{1}, 
\underline{\theta}_{2}) \hspace {8mm}
(i=0,1,2).
\end{equation*} 
Then, the above integrability condition is equivalent to $A=B=0$, where $A$ and $B$ 
are given by
\begin{align*}
A=&(f_{11})_{x_2}-(f_{12})_{x_1}+(f_{11})_{y}z_{2}+(f_{11})_{z_1}f_{12}+(f_{11})_{z_2}f_{22}
\\ &-(f_{12})_{y}z_{1}-(f_{12})_{z_1}f_{11}-(f_{12})_{z_2}f_{12}, \\
\\
B=&(f_{12})_{x_2}-(f_{22})_{x_1}+(f_{12})_{y}z_{2}+(f_{12})_{z_1}f_{12}+(f_{12})_{z_2}f_{22}
\\ &-(f_{22})_{y}z_{1}-(f_{22})_{z_1}f_{11}-(f_{22})_{z_2}f_{12}.
\end{align*}
\begin{remark}
Hereafter, we discuss only the second order PDE (1) with repect to $f_{ij}$ 
satisfying $A=B=0$. 
\end{remark}

A family of integral surfaces of $\mathcal{I}$ gives a 2-dimensional foliation on 
$J^{1}(\mathbb R^{2},\mathbb R)$. We describe an infinitesimal automorphism group of the 
foliation, and consider a principal bundle over $J^{1}(\mathbb R^{2},\mathbb R)$ 
with this group as a structure group.\par
The contact lift $\hat\phi$ of the scale transformation $\phi$ preserving 
$\mathcal{I}$ satisfies the following equations:
\begin{align}
\hat\phi^{*}\underline{\theta}_{0}&=a\underline{\theta}_{0} \hspace{30mm} (a\not=0),
\nonumber \\
\hat\phi^{*}\underline{\theta}_{1}&=b\underline{\theta}_{0}
+c\underline{\theta}_{1} \hspace{20mm} (c\not=0),  \nonumber \\
\hat\phi^{*}\underline{\theta}_{2}&=e\underline{\theta}_{0}
 +g\underline{\theta}_{2} \hspace{20mm} (g\not=0),  \\
\hat\phi^{*}\underline{\omega}_{1}&=h\underline{\omega}_{1} 
\hspace{29mm} (h\not=0), \nonumber \\
\hat\phi^{*}\underline{\omega}_{2}&=k\underline{\omega}_{2} 
\hspace{29mm} (k\not=0). \nonumber
\end{align}
The equation (4) can be written in the following form:
\begin{equation}
\begin{bmatrix}
\theta_{0} \\ \theta_{1} \\ \theta_{2} \\ \omega_{1} \\ \omega_{2}
\end{bmatrix}
=\begin{bmatrix}
a &0 &0 &0 &0 \\
b &c &0 &0 &0 \\
e &0 &g &0 &0 \\
0 &0 &0 &h &0 \\
0 &0 &0 &0 &k
\end{bmatrix}
\begin{bmatrix}
\underline{\theta}_{0} \\ \underline{\theta}_{1} \\
\underline{\theta}_{2} \\ \underline{\omega}_{1} \\ \underline{\omega}_{2}
\end{bmatrix}
\end{equation}
where, $a, b, c, e, g, h, k$ are functions. Thus we have linear transformations 
of coframes determined by $\hat\phi$.
Moreover, the lift $\hat\phi$ of the scale transformation satisfies:
\begin{align}
d\theta_{0} &\equiv -\theta_{1}\wedge\omega_{1}-\theta_{2}\wedge\omega_{2}   
\hspace{24mm}  (\text{mod } \theta_{0}), \nonumber \\ 
d\theta_{1} &\equiv 0 
\hspace{51mm} (\text{mod } \theta_{0},\theta_{1},\theta_{2}), \\
d\theta_{2} &\equiv 0 
\hspace{51mm} (\text{mod } \theta_{0},\theta_{1},\theta_{2}). \nonumber
\end{align}
These relations give conditions $a=ch=gk$. 
From these conditions, we get the linear transformations of coframes of the following form:
\begin{equation}
\begin{bmatrix}
\theta_{0} \\ \theta_{1} \\ \theta_{2} \\ \omega_{1} \\ \omega_{2}
\end{bmatrix}
=\begin{bmatrix}
ch &0 &0 &0 &0 \\
b &c &0 &0 &0 \\
e &0 &g &0 &0 \\
0 &0 &0 &h &0 \\
0 &0 &0 &0 &k
\end{bmatrix}
\begin{bmatrix}
\underline{\theta}_{0} \\ \underline{\theta}_{1} \\
\underline{\theta}_{2} \\ \underline{\omega}_{1} \\ \underline{\omega}_{2}
\end{bmatrix}.
\end{equation}
Therefore, 
we obtain the following 5-dimensional Lie group as infinitesimal automorphism group:
\begin{equation}
G:=\left.\left\{\begin{bmatrix}
ch &0 &0 &0 &0 \\
b &c &0 &0 &0 \\
e &0 &g &0 &0 \\
0 &0 &0 &h &0 \\
0 &0 &0 &0 &k
\end{bmatrix}
\in GL(5,\mathbb R) \right| ch=gk \right\}.
\end{equation}
Then, we choose the reduced $G$-bundle $\mathcal{F}_{G}$
 of the coframe bundle $\mathcal{F}_{GL}(\mathbb R^5)$ over $J^1(\mathbb R^2, \mathbb R)$. 
This bundle $\mathcal{F}_{G}$ is called $G$-structure associated with 
the second order PDE (1).
\section{Cartan's equivalence method}

\hspace{4mm} In the previous section, we introduced a $G$-structure $\mathcal{F}_{G}$ 
associated with the second order PDE (1). In this section we compute curvatures 
for the equivalence problem. For this purpose, 
we adopt Cartan's equivalence method ([3], [12], [18]).\par
First, we compute the structure equation on $\mathcal{F}_{G}$. From (7), we can choose\\ 
$(\theta_{0},\theta_{1},\theta_{2},\omega_{1},\omega_{2})$ as $\mathbb R^5$-valued 
tautological 1-form on $\mathcal{F}_{G}$. To obtain the structure equation, we 
compute the exterior derivative of the tautological 1-forms 
($\theta_{0},\theta_{1},\theta_{2},\omega_{1},\omega_{2}$).
\begin{equation}
\begin{split}
d\begin{bmatrix}
\theta_{0} \\ \theta_{1} \\ \theta_{2} \\ \omega_{1} \\ \omega_{2}
\end{bmatrix}
&=\begin{bmatrix}
\frac{dc}{c}+\frac{dh}{h} &0 &0 &0 &0 \\
\frac{db}{ch}-\frac{bdc}{c^{2}h} &\frac{dc}{c} &0 &0 &0 \\
\frac{de}{ch}-\frac{edg}{cgh} &0 &\frac{dg}{g} &0 &0 \\
0 &0 &0 &\frac{dh}{h} &0 \\
0 &0 &0 &0 &\frac{dk}{k}
\end{bmatrix}
\wedge\begin{bmatrix}
\theta_{0} \\ \theta_{1} \\ \theta_{2} \\ \omega_{1} \\ \omega_{2}
\end{bmatrix} \\
&+
\begin{bmatrix}
T_{1}\omega_{1}\wedge\theta_{0}+T_2\omega_{2}\wedge\theta_{0}
-\theta_{1}\wedge\omega_{1}-\theta_{2}\wedge\omega_{2} \\
\theta_{0}\wedge(T_3\omega_{1}+T_4\omega_{2})+\theta_{1}\wedge(T_5\omega_{1}
+T_6\omega_{2})+\theta_{2}\wedge(T_7\omega_{1}+T_8\omega_{2})  \\
\theta_{0}\wedge(T_9\omega_{1}+T_{10}\omega_{2})
+\theta_{1}\wedge(T_{11}\omega_{1}+T_{12}\omega_{2})
+\theta_{2}\wedge(T_{13}\omega_{1}+T_{14}\omega_{2}) \\0 \\0
\end{bmatrix},
\end{split}
\end{equation}

where,
\begin{align*}
 &T_1=-\frac{b}{ch},\hspace{1cm} T_2=-\frac{e}{ch} 
,\hspace{1cm} T_3=\frac{b^{2}}{(ch)^{2}}-\frac{(f_{11})_{y}}{h^2}+\frac{b(f_{11})_{z_1}}{ch^2}+\frac{e(f_{11})_{z_2}}{gh^2}, \\
&T_4=\frac{be}{(ch)^{2}}-\frac{(f_{12})_{y}}{hk}+\frac{b(f_{12})_{z_1}}{chk}
+\frac{e(f_{12})_{z_2}}{ch^2}, \hspace{1cm}
T_5=-\frac{b}{ch}-\frac{(f_{11})_{z_1}}{h},\\
&T_6=-\frac{(f_{12})_{z_1}}{k},\hspace{1cm}  T_7=-\frac{c(f_{11})_{z_2}}{gh}, \hspace{1cm}
T_8=-\frac{b}{ch}-\frac{(f_{12})_{z_2}}{h}, \\
&T_9=\frac{be}{(ch)^{2}}-\frac{g(f_{12})_{y}}{ch^2}+\frac{bg(f_{12})_{z_1}}{(ch)^2}
+\frac{e(f_{12})_{z_2}}{ch^2},\\ 
&T_{10}=\frac{e^2}{(ch)^{2}}-\frac{g(f_{22})_{y}}{chk}+\frac{bg(f_{22})_{z_1}}{c^{2}hk}
+\frac{e(f_{22})_{z_2}}{chk}, \\
&T_{11}=-\frac{e}{ch}-\frac{g(f_{12})_{z_1}}{ch},\hspace{1cm}
T_{12}=-\frac{g(f_{22})_{z_1}}{ck},\hspace{1cm} T_{13}=-\frac{(f_{12})_{z_2}}{h},\\
&T_{14}=-\frac{e}{ch}-\frac{(f_{22})_{z_2}}{k}.
\end{align*}
\begin{remark}
We put $\omega$:=($\theta_{0},\theta_{1},\theta_{2},\omega_{1},\omega_{2}$) and write the 
structure equation (9):
\begin{center}
$d\omega=-\theta\wedge\omega+T\omega\wedge\omega.$
\end{center}
In the above, we note that $\theta$ is a $\mathfrak g$-valued 1-form and $T\omega\wedge\omega$ is a $\mathbb R^{5}$-valued 2-form. In fact, 
\begin{center}
$d\omega=d(g\underline\omega)=dg\cdot g^{-1}\wedge\omega+T\omega\wedge\omega,$
\end{center}
where $g \in G$ and 
$\underline\omega=(\underline{\theta_{0}}, \underline{\theta_{1}}, 
\underline{\theta_{2}}, \underline{\omega_{1}}, \underline{\omega_{2}})$.
The above equation shows that $\theta$ is the Maurer-Cartan form. In the structure equation (9),
 each component of $\theta$ is called the pseudo-connection form and 
$T\omega\wedge\omega$ is called the torsion 2-form, and coefficient functions of 2-forms in each component of $T\omega\wedge\omega$ are called torsions ([6]).
\end{remark}
To simplify the structure equation (9), we set:
\begin{align*}
\alpha:=&\frac{dc}{c}-\frac{b}{ch}\omega_{1}-\frac{e}{ch}\omega_{2}, \\
\beta:=&\frac{db}{ch}-\frac{bdc}{c^{2}h}-
\left\{\frac{b^2}{(ch)^2}-\frac{(f_{11})_{y}}{h^2}+\frac{b(f_{11})_{z_1}}{ch^2}
+\frac{e(f_{11})_{z_2}}{gh^2}\right\}\omega_{1} \\
&-\left\{\frac{be}{(ch)^2}-\frac{(f_{12})_{y}}{hk}+\frac{b(f_{12})_{z_1}}{chk}
+\frac{e(f_{12})_{z_2}}{ch^2}\right\}\omega_{2}, \\
\varepsilon:=&\frac{de}{ch}-\frac{edg}{cgh}-
\left\{\frac{be}{(ch)^2}-\frac{g(f_{12})_{y}}{ch^2}+\frac{bg(f_{12})_{z_1}}{(ch)^2}
+\frac{e(f_{12})_{z_2}}{ch^2}\right\}\omega_{1} \\
&-\left\{\frac{e^2}{(ch)^2}-\frac{g(f_{22})_{y}}{chk}+\frac{bg(f_{22})_{z_1}}{c^2hk}
+\frac{e(f_{22})_{z_2}}{chk}\right\}\omega_{2}, \\
\delta:=&\frac{dg}{g}-\frac{b}{ch}\omega_{1}-\frac{e}{ch}\omega_{2},\quad
\gamma:=\frac{dh}{h}, \quad
\psi:=\frac{dk}{k}.
\end{align*}

By substituting the above terms into the equation (9), we get the following proposition.
\begin{proposition}
The structure equation on $\mathcal {F}_{G}$ is written as:
\begin{equation}
\begin{split}
d\begin{bmatrix}
\theta_{0} \\ \theta_{1} \\ \theta_{2} \\ \omega_{1} \\ \omega_{2}
\end{bmatrix}
&=\begin{bmatrix}
\alpha+\gamma & 0& 0& 0& 0 \\
\beta & \alpha& 0& 0 &0 \\
\varepsilon &0 &\delta &0 & 0 \\
0& 0& 0& \gamma &0 \\
0& 0& 0& 0& \psi
\end{bmatrix}
\wedge\begin{bmatrix}
\theta_{0} \\ \theta_{1} \\ \theta_{2} \\ \omega_{1} \\ \omega_{2}
\end{bmatrix} \\
&+\begin{bmatrix}
-\theta_{1}\wedge\omega_{1}-\theta_{2}\wedge\omega_{2} \\
L_1\theta_{1}\wedge\omega_{1}+L_2\theta_{1}\wedge\omega_{2}+L_3\theta_{2}\wedge\omega_{1}
+L_4\theta_{2}\wedge\omega_{2}  \\
L_2\theta_{1}\wedge\omega_{1}+L_5\theta_{1}\wedge\omega_{2}+L_4\theta_{2}\wedge\omega_{1}
+L_6\theta_{2}\wedge\omega_{2} \\
0 \\0
\end{bmatrix},
\end{split}
\end{equation}

where, 
\begin{align*}
&L_1:=-\frac{2b}{ch}-\frac{(f_{11})_{z_1}}{h},\hspace{1cm}
L_2:=-\frac{e}{ch}-\frac{(f_{12})_{z_1}}{k},\hspace{1cm}
L_3:=-\frac{c(f_{11})_{z_2}}{gh}, \\
&L_4:=-\frac{b}{ch}-\frac{(f_{12})_{z_2}}{h},\hspace{1cm}
L_5:=-\frac{g(f_{22})_{z_1}}{ck},\hspace{1cm}
L_6:=-\frac{2e}{ch}-\frac{(f_{22})_{z_2}}{k},\\
&  \alpha+\gamma=\delta+\psi.
\end{align*}
\end{proposition}
\begin{remark}
In (10), some torsions in the structure equation (9) are absorved.
\end{remark}
To eliminate the ambiguity of the pseudo-connection forms, we need to choose a reduction 
of $G$-structure $\mathcal {F}_{G}$. 
We choose the reduction of $G$-structure by setting $L_2=L_4=0$. 
We denote this reduced bundle by $\mathcal{F}_{G_1}$, where $G_1$ is the following 
3-dimensional Lie group:
\begin{equation*}
G_1:=\left.\left\{\begin{bmatrix}
ch &0 &0 &0 &0 \\
0 &c &0 &0 &0 \\
0 &0 &g &0 &0 \\
0 &0 &0 &h &0 \\
0 &0 &0 &0 &k
\end{bmatrix}
\in GL(5,\mathbb R) \right| ch=gk \right\}.
\end{equation*}
We have the tautological 1-form on $\mathcal{F}_{G_1}$ given by:
\begin{equation*}
\begin{bmatrix}
\hat\theta_{0} \\ \hat\theta_{1} \\ \hat\theta_{2} \\ \hat\omega_{1} \\ \hat\omega_{2}
\end{bmatrix}
=\begin{bmatrix}
ch\underline{\theta_{0}} \\ 
-c(f_{12})_{z_2}\underline{\theta_{0}}+c\underline{\theta_{1}} \\
-g(f_{12})_{z_1}\underline{\theta_{0}}+g\underline{\theta_{2}} \\
 h\underline{\omega_{1}} \\ k\underline{\omega_{2}}
\end{bmatrix}.
\end{equation*}

Then, the structure equation on $\mathcal{F}_{G_1}$ is given by
\begin{equation*}
\begin{split}
\hspace{-0.5cm} 
&d\begin{bmatrix}
\hat\theta_{0} \\ \hat\theta_{1} \\ \hat\theta_{2} \\ \hat\omega_{1} \\ \hat\omega_{2}
\end{bmatrix}
=\begin{bmatrix}
\alpha+\gamma & 0& 0& 0& 0 \\
0 & \alpha& 0& 0 &0 \\
0 &0 &\delta &0 & 0 \\
0& 0& 0& \gamma &0 \\
0& 0& 0& 0& \psi
\end{bmatrix}
\wedge\begin{bmatrix}
\hat\theta_{0} \\ \hat\theta_{1} \\ \hat\theta_{2} \\ \hat\omega_{1} \\ \hat\omega_{2}
\end{bmatrix} \\
&+\resizebox{0.95\hsize}{!}{$
\begin{bmatrix}
M_{12}\hat\omega_{1}\wedge\hat\theta_{0}+M_{11}\hat\omega_{2}\wedge\hat\theta_{0}
-\hat\theta_{1}\wedge\hat\omega_{1}-\hat\theta_{2}\wedge\hat\omega_{2} \\
M_1\hat\theta_{2}\wedge\hat\omega_{1}+M_2\hat\theta_{1}\wedge\hat\theta_{0}
+M_3\hat\theta_{2}\wedge\hat\theta_{0}
+M_4\hat\omega_{1}\wedge\hat\theta_{0}+M_5\hat\omega_{2}\wedge\hat\theta_{0}
+M_{10}\hat\omega_{1}\wedge\hat\theta_{1}+M_{11}\hat\omega_{2}\wedge\hat\theta_{1} \\
M_6\hat\theta_{1}\wedge\hat\omega_{2}+M_7\hat\theta_{1}\wedge\hat\theta_{0}
+M_2\hat\theta_{2}\wedge\hat\theta_{0}+M_8\hat\omega_{1}\wedge\hat\theta_{0}
+M_9\hat\omega_{2}\wedge\hat\theta_{0}+M_{12}\hat\omega_{1}\wedge\hat\theta_{2}
+M_{13}\hat\omega_{2}\wedge\hat\theta_{2} \\
0 \\0
\end{bmatrix}
$}
\end{split}
\end{equation*}
where, 
\begin{align*}
\alpha&=\frac{dc}{c} ,\hspace{1cm}
\delta=\frac{dg}{g} ,\hspace{1cm}
\gamma=\frac{dh}{h},\hspace{1cm} \psi=\frac{dk}{k},
\end{align*}
\[\hspace{-1cm}
M_1=-\frac{c(f_{11})_{z_2}}{gh}, \hspace{1cm} M_2=-\frac{(f_{12})_{z_2z_1}}{ch}
,\hspace{1cm} M_3=-\frac{(f_{12})_{z_2z_2}}{gh}, 
\]
\begin{multline*}
M_4=-\frac{1}{h^2} \bigl\{ (f_{12})_{z_2}^2-(f_{11})_y
-(f_{12})_{z_2}(f_{11})_{z_1}-(f_{11})_{z_2}(f_{12})_{z_1}\\
+(f_{12})_{z_{2}x_{1}}+(f_{12})_{z_{2}y}z_{1}+(f_{12})_{z_{2}z_{1}}f_{11}
+(f_{12})_{z_{2}z_{2}}f_{21} \bigr\}, 
\end{multline*}
\begin{multline*}
M_5=\frac{1}{hk} \bigl\{(f_{12})_{y}+(f_{12})_{z_{2}}(f_{12})_{z_{1}}-(f_{12})_{z_{2}x_{2}}
-(f_{12})_{z_2y}z_2-(f_{12})_{z_{2}z_{1}}f_{12}-(f_{12})_{z_{2}z_{2}}f_{22}\bigr\}, 
\end{multline*}
\[\hspace{-8cm} M_6=-\frac{g(f_{22})_{z_1}}{ck} ,\hspace{1cm} M_7=-\frac{(f_{12})_{z_1z_1}}{ck}, \]
\begin{multline*}
M_8=\frac{1}{hk} \bigl\{(f_{12})_{y}+(f_{12})_{z_{1}}(f_{12})_{z_{2}}-(f_{12})_{z_{1}x_{1}}
-(f_{12})_{z_1y}z_1-(f_{12})_{z_{1}z_{1}}f_{11}-(f_{12})_{z_{1}z_{2}}f_{21}\bigr\}, 
\end{multline*}
\begin{multline*}
M_9=-\frac{1}{k^2} \bigl\{ (f_{12})_{z_1}^2-(f_{22})_y
-(f_{12})_{z_2}(f_{22})_{z_1}-(f_{12})_{z_1}(f_{22})_{z_2}\\
\hspace{-5cm} +(f_{12})_{z_{1}x_{2}}+(f_{12})_{z_{1}y}z_{2}+(f_{12})_{z_{1}z_{1}}f_{12}
+(f_{12})_{z_{1}z_{2}}f_{22} \bigr\}. 
\end{multline*}
\begin{align*}
M_{10}=\frac{1}{h}\left\{(f_{11})_{z_1}-(f_{12})_{z_2}\right\}, 
\hspace{1cm} M_{11}=\frac{(f_{12})_{z_1}}{k},\\
M_{12}=\frac{(f_{12})_{z_2}}{h}, \hspace{1cm} 
M_{13}=\frac{1}{k}\left\{(f_{22})_{z_2}-(f_{12})_{z_1}\right\}.
\end{align*}
By absorption of torsions $M_2,\ M_{10},\ M_{11},\ M_{12},\ M_{13}$, 
we obtain the following:
\begin{proposition}
We have the following structure equation on $\mathcal{F}_{G_1}$.
\begin{equation}
\begin{split}
d\begin{bmatrix}
\hat\theta_{0} \\ \hat\theta_{1} \\ \hat\theta_{2} \\ \hat\omega_{1} \\ \hat\omega_{2}
\end{bmatrix}
&=
\begin{bmatrix}
\hat\alpha+\hat\gamma & 0& 0& 0& 0 \\
0 & \hat\alpha& 0& 0 &0 \\
0 &0 &\hat\delta &0 & 0 \\
0& 0& 0& \hat\gamma &0 \\
0& 0& 0& 0& \hat\psi
\end{bmatrix}
\wedge\begin{bmatrix}
\hat\theta_{0} \\ \hat\theta_{1} \\ \hat\theta_{2} \\ \hat\omega_{1} \\ \hat\omega_{2}
\end{bmatrix} \\
&+
\begin{bmatrix}
\hat\omega_{1}\wedge\hat\theta_{1}+\hat\omega_{2}\wedge\hat\theta_{2} \\
M_1\hat\theta_{2}\wedge\hat\omega_{1}+M_3\hat\theta_{2}\wedge\hat\theta_{0}
+M_4\hat\omega_{1}\wedge\hat\theta_{0}+M_5\hat\omega_{2}\wedge\hat\theta_{0} \\
M_6\hat\theta_{1}\wedge\hat\omega_{2}+M_7\hat\theta_{1}\wedge\hat\theta_{0}
+M_8\hat\omega_{1}\wedge\hat\theta_{0}+M_9\hat\omega_{2}\wedge\hat\theta_{0} \\
0 \\0
\end{bmatrix}
\end{split}
\end{equation}
where we set
\begin{align*}
\hat\alpha&:=\alpha-M_{2}\hat\theta_{0}+M_{10}\hat\omega_{1}+M_{11}\hat\omega_{2}, \\
\hat\gamma&:=\gamma+(M_{12}-M_{10})\hat\omega_{1}, \\
\hat\delta&:=\delta-M_{2}\hat\theta_{0}+M_{12}\hat\omega_{1}+M_{13}\hat\omega_{2}, \\
\hat\psi&:=\psi+(M_{11}-M_{13})\hat\omega_{2}.
\end{align*}
\end{proposition}
We note that the structure equation (11) defines uniquely the pseudo-connection forms 
$\hat\alpha$, $\hat\gamma$, $\hat\delta$, $\hat\psi$. 
Hence, we can obtain the invariant 1-forms
$(\hat\theta_0,\hat\theta_{1},\hat\theta_{2},\hat\omega_{1},\hat\omega_{2},
\hat\alpha,\hat\gamma,\hat\psi)$ 
on $\mathcal{F}_{G_1}$. To consider the curvatures for the equivalence problem,
 we need to use the $\left\{e\right\}$-structure by 
choosing a prolongation of $\mathcal{F}_{G_1}$. 
Then, we obtain the following structure equation on the $\left\{e\right\}$-structure 
by taking the exterior derivation of tautological 1-forms
$(\hat\theta_0,\hat\theta_{1},\hat\theta_{2},\hat\omega_{1},\hat\omega_{2},
\hat\alpha, \hat\gamma,\hat\psi)$:
\begin{equation*}
\resizebox{0.95\hsize}{!}{$
d\begin{bmatrix}
\hat\theta_0 \\ \hat\theta_{1} \\ \hat\theta_{2} \\ \hat\omega_{1} \\ \hat\omega_{2}
\\   \hat\alpha \\ \hat\gamma \\ \hat\psi
\end{bmatrix}
=\begin{bmatrix}
(\hat\alpha+\hat\gamma)\wedge\hat\theta_0+\hat\omega_{1}\wedge\hat\theta_{1}
+\hat\omega_{2}\wedge\hat\theta_{2} \\
\hat\alpha\wedge\hat\theta_1+M_1\hat\theta_{2}\wedge\hat\omega_{1}
+M_3\hat\theta_{2}\wedge\hat\theta_{0}+M_4\hat\omega_{1}\wedge\hat\theta_{0}
+M_5\hat\omega_{2}\wedge\hat\theta_{0} \\
(\hat\alpha+\hat\gamma-\hat\psi)\wedge\hat\theta_2+M_6\hat\theta_{1}\wedge\hat\omega_{2}
+M_7\hat\theta_{1}\wedge\hat\theta_{0}+M_8\hat\omega_{1}\wedge\hat\theta_{0}
+M_9\hat\omega_{2}\wedge\hat\theta_{0} \\
\hat\gamma\wedge\hat\omega_1 \\
 \hat\psi\wedge\hat\omega_2 \\
S_1\hat\omega_{1}\wedge\hat\theta_{0}+S_2\hat\omega_{2}\wedge\hat\theta_{0}
+S_3\hat\theta_{1}\wedge\hat\theta_{0}+S_4\hat\theta_{2}\wedge\hat\theta_{0}
+S_5\hat\omega_{1}\wedge\hat\theta_{1}+S_6\hat\omega_{1}\wedge\hat\omega_{2}
+S_7\hat\theta_{2}\wedge\hat\omega_{1}-M_7\hat\theta_{1}\wedge\hat\omega_{2} \\
S_8\hat\omega_{1}\wedge\hat\omega_{2}+S_9\hat\omega_{1}\wedge\hat\theta_{0}
+S_5\hat\theta_{1}\wedge\hat\omega_{1}+S_{10}\hat\theta_{2}\wedge\hat\omega_{1} \\
S_{11}\hat\omega_{1}\wedge\hat\omega_{2}+S_{12}\hat\omega_{2}\wedge\hat\theta_{0} 
+S_{13}\hat\theta_{1}\wedge\hat\omega_{2}+S_{14}\hat\theta_{2}\wedge\hat\omega_{2}
\end{bmatrix}
$}.
\end{equation*}
Here, the torsions $M_i$ are given by previous page. To write down the torsions explicitly, we use the dual frame of the coframe 
$(\underline\theta_{0},\ \underline\theta_{1},\ \underline\theta_{2},\ \underline\omega_{1}
,\ \underline\omega_{2})$:
\begin{align*}
\partial_{\underline\theta_{0}}:&=\frac{\partial}{\partial y},\quad
\partial_{\underline\theta_{1}}:=\frac{\partial}{\partial z_1},\quad
\partial_{\underline\theta_{2}}:=\frac{\partial}{\partial z_2},\\
\partial_{\underline\omega_{1}}:&=\frac{\partial}{\partial x_1}
+z_1\frac{\partial}{\partial y}+f_{11}\frac{\partial}{\partial z_1}
+f_{12}\frac{\partial}{\partial z_2},\\
\partial_{\underline\omega_{2}}:&=\frac{\partial}{\partial x_2}
+z_2\frac{\partial}{\partial y}+f_{21}\frac{\partial}{\partial z_1}
+f_{22}\frac{\partial}{\partial z_2}.
\end{align*}
By using the Frobenius condition $A=B=0$ 
(i.e. $(f_{11})_{\underline\omega_2}=(f_{12})_{\underline\omega_1},\ 
(f_{12})_{\underline\omega_2}=(f_{22})_{\underline\omega_1}$), each torsions of the above structure equation can be written
as follows.
\begin{align*}
M_{1}&=-\frac{c}{gh}(f_{11})_{\underline\theta_{2}}, \quad \quad \quad
M_{3}=-\frac{1}{gh}(f_{12})_{\underline\theta_{2}\underline\theta_{2}}, \\
M_{4}&=-\frac{1}{h^2}\left\{(f_{11})_{\underline\theta_2\underline\omega_2}-2(f_{11})_{\underline\theta_2}(f_{12})_{\underline\theta_1}
+(f_{11})_{\underline\theta_2}(f_{22})_{\underline\theta_2}\right\}, \\
M_5&=\frac{1}{hk}\left\{(f_{12})_{\underline\theta_0}
+(f_{12})_{\underline\theta_2}(f_{12})_{\underline\theta_1}
-(f_{12})_{\underline\theta_2\underline\omega_2}\right\},\\
M_6&=-\frac{g}{ck}(f_{22})_{\underline\theta_1}, \quad \quad \quad
M_7=-\frac{1}{ck}(f_{12})_{\underline\theta_1\underline\theta_1}, \\
M_8&=\frac{1}{hk}\left\{(f_{12})_{\underline\theta_0}
+(f_{12})_{\underline\theta_1}(f_{12})_{\underline\theta_2}
-(f_{12})_{\underline\theta_1\underline\omega_1}\right\}, \\
M_9&=-\frac{1}{k^2}\left\{-2(f_{12})_{\underline\theta_2}(f_{22})_{\underline\theta_1}
+(f_{22})_{\underline\theta_1\underline\omega_1}
+(f_{11})_{\underline\theta_1}(f_{22})_{\underline\theta_1}\right\}, \\
S_1&=\frac{1}{ch^2}\bigl\{(f_{11})_{\underline\theta_2\underline\theta_1
\underline\omega_2}+(f_{11})_{\underline\theta_2\underline\theta_2}
(f_{22})_{\underline\theta_1}+(f_{11})_{\underline\theta_2\underline\theta_1}
(f_{22})_{\underline\theta_2}-(f_{12})_{\underline\theta_2\underline\theta_1}
(f_{11})_{\underline\theta_2}\\ &-(f_{12})_{\underline\theta_2\underline\theta_2} 
(f_{12})_{\underline\theta_2}-(f_{12})_{\underline\theta_2\underline\theta_1}
(f_{11})_{\underline\theta_1}+2(f_{12})_{\underline\theta_2\underline\theta_1}
(f_{12})_{\underline\theta_2}
\bigr\}, \\
S_2&=\frac{1}{chk}\bigl\{(f_{12})_{\underline\theta_2\underline\theta_1\underline\omega_2}
-(f_{12})_{\underline\theta_1\underline\theta_0}
-(f_{12})_{\underline\theta_1\underline\theta_1}(f_{12})_{\underline\theta_2}\bigr\}, \\
 S_3&=\frac{(f_{12})_{\underline\theta_2\underline\theta_1\underline\theta_1}}{c^2h},
 \hspace{1.5cm}
S_4=\frac{(f_{12})_{\underline\theta_2\underline\theta_1\underline\theta_2}}{cgh},
 \hspace{1.5cm}
S_5=\frac{2(f_{12})_{\underline\theta_2\underline\theta_1}
-(f_{11})_{\underline\theta_1\underline\theta_1}}{ch},\\
\hspace{-3cm} S_6&=\frac{1}{hk}\left\{-(f_{12})_{\underline\theta_0}
-(f_{12})_{\underline\theta_1}(f_{12})_{\underline\theta_2}
+(f_{11})_{\underline\theta_2}(f_{22})_{\underline\theta_1}
+(f_{12})_{\underline\theta_2\underline\omega_2} \right\},\\
S_7&=\frac{(f_{11})_{\underline\theta_1\underline\theta_2}
-(f_{12})_{\underline\theta_2\underline\theta_2}}{gh}, \hspace{2cm}
S_8=\frac{1}{hk}\left\{(f_{11})_{\underline\theta_1\underline\omega_2}
-2(f_{12})_{\underline\theta_2\underline\omega_2} \right\},\\
S_9&=\frac{1}{ch^2}\bigl\{(f_{11})_{\underline\theta_1\underline\theta_0}
-2(f_{12})_{\underline\theta_2\underline\theta_0}
+(f_{11})_{\underline\theta_1\underline\theta_1}(f_{12})_{\underline\theta_2}\\
&\quad\quad+(f_{11})_{\underline\theta_1\underline\theta_2}(f_{12})_{\underline\theta_1}
-2(f_{12})_{\underline\theta_1\underline\theta_2}(f_{12})_{\underline\theta_2}
-2(f_{12})_{\underline\theta_2\underline\theta_2}(f_{12})_{\underline\theta_1} \bigr\},\\
S_{10}&=\frac{-(f_{11})_{\underline\theta_1\underline\theta_2}
+2(f_{12})_{\underline\theta_2\underline\theta_2}}{gh},\quad\quad 
S_{11}=\frac{1}{hk}\bigl\{2(f_{12})_{\underline\theta_1\underline\omega_1}
-(f_{22})_{\underline\theta_2\underline\omega_1} \bigr\},\\
S_{12}&=\frac{1}{chk}\bigl\{-2(f_{12})_{\underline\theta_1\underline\theta_0}
-2(f_{12})_{\underline\theta_1\underline\theta_1}(f_{12})_{\underline\theta_2}
-2(f_{12})_{\underline\theta_1\underline\theta_2}(f_{12})_{\underline\theta_1}
+(f_{22})_{\underline\theta_2\underline\theta_0}\\
&\quad\quad+(f_{22})_{\underline\theta_1\underline\theta_2}(f_{12})_{\underline\theta_2}
+(f_{22})_{\underline\theta_2\underline\theta_2}(f_{12})_{\underline\theta_1} \bigr\},\\
S_{13}&=\frac{2(f_{12})_{\underline\theta_1\underline\theta_1}
-(f_{22})_{\underline\theta_1\underline\theta_2}}{ck},\quad\quad 
S_{14}=\frac{2(f_{12})_{\underline\theta_1\underline\theta_2}
-(f_{22})_{\underline\theta_2\underline\theta_2}}{gk}.
\end{align*}
In the above torsions, there are the following relations.
\begin{proposition}
Torsions $M_4,\ M_9,\ S_3,\ S_4,\ S_7,\ S_{10},\ S_{13}$ are given by;
\begin{align*}
M_4&=-\frac{1}{h^2}\left\{\frac{-gh}{c}(M_1)_{\underline\omega_2}
+\frac{2gh}{c}M_1(f_{12})_{\underline\theta_1}
-\frac{gh}{c}M_1(f_{22})_{\underline\theta_2}\right\}, \\
M_9&=-\frac{1}{k^2}\left\{-\frac{ck}{g}(M_6)_{\underline\omega_1}
-\frac{ck}{g}M_6(f_{11})_{\underline\theta_1}
+\frac{2ck}{g}M_6(f_{12})_{\underline\theta_2}\right\},\\
S_3&=-\frac{k}{ch}(M_7)_{\underline\theta_2}, \hspace{1cm}
S_4=-\frac{1}{c}(M_3)_{\underline\theta_1}, \hspace{1cm}
S_7=-\frac{1}{c}(M_1)_{\underline\theta_1}+M_3, \\
S_{10}&=-\frac{1}{c}(M_1)_{\underline\theta_1}+2M_3, \hspace{1cm}
S_{13}=-2M_7+\frac{1}{g}(M_6)_{\underline\theta_2}.
\end{align*}
\end{proposition}
Hence, the vanishing of $M_4,\ M_9,\ S_3,\ S_4,\ S_7,\ S_{10},\ S_{13}$ is 
given by vanishing of other curvatures. 
By the theory of $G$-structure ([12], [18]), a vanishing condition of curvatures
$M_{i},\ S_{j}$ ($i=1,3,5,6,7,8,\ j=1,2,5,6,8,9,11,12,14$) gives the 
following theorem.
\begin{theorem}
Suppose that the second order PDE {\rm(1)} satisfies the integrability 
condition $A=B=0$. Then, the equation {\rm(1)}
 is {\rm(}locally{\rm)} equivalent to the flat equation 
under lifts of scale transformations if and only if curvatures $M_i,\ S_j$ vanish.
\end{theorem}

First, it is easy to check that the functions $f_{ij}$ satisfying $A=B=M_{i}=S_{j}=0$
 are written as quadratic polynomials in $z_1,\ z_2$.
Hence, if there is a polynomial $z_1,\ z_2$ of degree three among $f_{ij}$, then corresponding 
equation (1) is not equivalent to the flat equation under lifts of scale 
transformations.\par
Next, we give some examples of equation which is equivalent to the flat equation.
To show the vanishing condition of the curvatures more explicitly, we consider the functions $f_{ij}$ given by:
\begin{center}
$f_{11}=P(x_1,x_2,y),\hspace{1cm} f_{12}=Q(x_1,x_2,y),\hspace{1cm} f_{22}=R(x_1,x_2,y).$
\end{center}
Then, Theorem 3.6 gives the following Corollary.
\begin{corollary}
Suppose that the functions $f_{ij}$ in $(1)$ are given in the above form.
Then the equation $(1)$ is {\rm(}locally{\rm)} equivalent to the flat equation under 
the lifts of scale transformations if and only if $P_{y}=Q_{y}=R_{y}=0$, 
 $P_{x_2}=Q_{x_1}$, $Q_{x_2}=R_{x_1}$.
\end{corollary}
\begin{remark}
The conditions $P_{y}=Q_{y}=R_{y}=0$, 
 $P_{x_2}=Q_{x_1}$, $Q_{x_2}=R_{x_1}$ in Corollary 3.7 are obtained by the 
integrability condition $A=B=0$. 
Namely, a vanishing condition of curvatures (i.e. $M_{i}=S_{j}=0$) is absorved into the
integrability condition. 
Therefore, it is shown that the second order PDE (1) for the functions $f_{ij}$ given by
 the above form are equivalent to the flat equations if and only if it is integrable.
\end{remark}

\section{Duality associated with differential equations}
\hspace{4mm}
In this section, we discuss a duality between the coordinate space and the solution space 
associated with the following flat equation; 
\begin{equation}\label{flat equation}
\frac{\partial ^{2} y}{\partial x_i \partial x_j}=0 \hspace{4mm} (1\leq i,j \leq 2). 
\end{equation}
For the purpose, we consider the following double fibration.
\begin{equation}\label{flat model}
\setlength{\unitlength}{0.7mm}
\begin{picture}(100,50)(-30,0)
\put(0,0){\makebox(-5,10)[r]{$\mathbb R^3:=\left\{(x_1, x_2, y)\right\}$}}
\put(16,30){\makebox(-10,10){$J^1(\mathbb R^2,\mathbb R)$}}
\put(35,0){\makebox(10,10)[l]{$\mathbb R^3:=\left\{(a,b,c)\right\}$}}
\put(5,17){\makebox(-40,10){$\pi_1$}}
\put(33,16){\makebox(0,10){$\pi_2$}}
\put(1,30){\vector(-1,-1){20}}
\put(20,30){\vector(1,-1){20}}
\end{picture}
\end{equation}
where, projections $\pi_{1},\ \pi_{2}$ are defined by 
\begin{align*}
\pi_{1}(x_1,x_2,y,z_1,z_2)&=(x_1,x_2,y), \\
\pi_{2}(x_1,x_2,y,z_1,z_2)&=(z_1,z_2,y-z_{1}x_{1}-z_{2}x_{2}).
\end{align*}
We call the double fibration (\ref{flat model}) 
the model space of the flat equation or flat model space.
In this fibration, we regard the left base space as a coordinate space 
$\mathbb R^3:=\left\{(x_1,x_2,y)\right\}$, and a right base space as a 
solution space $\mathbb R^3:=\left\{(a,b,c)\right\}$.
Solutions of (\ref{flat equation}) are written as $y=ax_1+bx_2+c$ 
for real parameters $a,b,c$. 
Graphs of solutions are planes on $\mathbb R^3$ or $J^1(\mathbb R^2,\mathbb R)$, 
and the 3-parameter family of solutions yields 
a 2-dimensional foliation on $J^1(\mathbb R^2,\mathbb R)$. 
Then the leaf space of this foliation is interpreted as a solution space of 
 (\ref{flat equation}).
We discuss the compactification of the flat model space. The 
fibration (\ref{flat model}) can be embedded naturally into the 
following (global) double fibration:
\setlength{\unitlength}{0.7mm}
\begin{equation}\label{flag}
\begin{picture}(100,50)(-30,0)
\put(-40,0){\makebox(20,10)[r]{$\mathbb P(V)$}}
\put(0,30){\makebox(20,10){$\mathbb F_{V}(1,3)$}}
\put(20,0){\makebox(20,10)[l]{$\mathbb P(V^*)\cong Gr(3, 4)$}}
\put(-25,17){\makebox(5,10){$\pi_1$}}
\put(30,16){\makebox(35,10){$\pi_2$}}
\put(1,30){\vector(-1,-1){20}}
\put(20,30){\vector(1,-1){20}}
\end{picture}
\end{equation}
where, $V=\mathbb R^4$, and 
$Gr(3, 4)$ is a Grassmannian manifold and $\mathbb F_{V}(1, 3)$ is a flag variety: 
\begin{align*}             
&Gr(3,4)=\left\{E \ |\ E {\rm \ is\ a\ hyperplane\ of}\ V:=\mathbb R^4 \right\}, \\
&\mathbb F_{V}(1,3)=\left\{(l, E)\ |\ l \in \mathbb RP^3,\ E \in Gr(3, 4)\cong \mathbb RP^3 ,\ l\subset E \right\}.
\end{align*}
The, projections $\pi_1, \pi_2$ are defined by 
\begin{center}
$\pi_1([u], H):=[u], \quad\quad
\pi_2([u], H):=[f_H],$
\end{center}
where, $f_H$ is a linear functional satisfying $ker(f_H)=H$ of 
$V^*\backslash \left\{0\right\}$.
(Since $f_H$ is uniquely defined up to scalar multiplication, $\pi_2$ is well-defined.)
The double fibration (14) 
 do not depend on coordinate transformation group $\mathcal G$. 
So, we introduce the flat model space depending on $\mathcal G$ .\par
 We fix a coordinate transformation group $\mathcal G\subset {\rm Diff}(\mathbb R^3)$. 
First, we define the following symmetry group. $([11],\ [12])$
\begin{definition}
Let $G$ be an isotropy subgroup of the flat equation (\ref{flat equation}) in $\mathcal G$.
 This group $G$ is called symmetry group of the flat equation for $\mathcal G$.
\end{definition}
In the case of $\mathcal G={\rm Diff}(\mathbb R^3)$, the symmetry group is $SL(4, \mathbb R)$ 
and the action on the coordinate space $\mathbb R^3$ is given by:
\vspace{0.2cm}

For 
$\begin{pmatrix}
a_1&a_2&a_3&a_4 \\
b_1&b_2&b_3&b_4 \\
c_1&c_2&c_3&c_4 \\
d_1&d_2&d_3&d_4 \\
\end{pmatrix}
\in SL(4, \mathbb R)$, 
\begin{equation}
(x_1, x_2, y) \mapsto \left
(\frac{a_1x_1+a_2x_2+a_3y+a_4}{d_1x_1+d_2x_2+d_3y+d_4}, 
\frac{b_1x_1+b_2x_2+b_3y+b_4}{d_1x_1+d_2x_2+d_3y+d_4}, 
\frac{c_1x_1+c_2x_2+c_3y+c_4}{d_1x_1+d_2x_2+d_3y+d_4} \right).
\end{equation}
Next, we introduce subgroups of $SL(4, \mathbb R)$ as follows:
\begin{align*}
H_i:&=\left\{g \in SL(4, \mathbb R)\ |\ g[e_{i}]=[e_{i}] \right\}, \\
\overline H_i:&=\left\{g \in SL(4, \mathbb R)\ |\ ^{t}g^{-1}[e_{i}]=[e_{i}] \right\},
\end{align*}
where, $e_i$ (i=1,$\cdot\cdot\cdot\cdot$,4) are standard basis of $\mathbb R^4$, and 
$[e_{i}]$ are corresponding elements in $RP^3$.
The subgroups $H_i$ are isotropy subgroups which preserve lines $[e_i]$, and
the subgroups $\overline H_i$ are isotropy subgroups which preserve hyperplanes
spanned by $e_j\ (j\not=i)$ respectively. We used Cartan involution 
$\tilde\theta(g)=^{t}g^{-1}$ in the definition of $\overline H_i$. 
We consider the following double fibration.
\begin{equation}
\setlength{\unitlength}{0.7mm}
\begin{picture}(100,50)(-30,0)
\put(3,0){\makebox(-20,10)[r]{$G/(G\cap H_4$)}}
\put(16,30){\makebox(-10,10){$G/(G\cap H)$}}
\put(39,0){\makebox(20,10)[l]{$G/(G\cap \overline H_1)$}}
\put(0,17){\makebox(-75,10){$(G\cap H_4)/(G\cap H)$}}
\put(40,16){\makebox(40,10){$(G\cap \overline H_1)/(G\cap H)$}}
\put(1,30){\vector(-1,-1){20}}
\put(20,30){\vector(1,-1){20}}
\end{picture}
\end{equation}
where, $H=H_4\cap \overline H_1$. 
We call this fibration as a model space of the flat equation (\ref{flat equation})
 with respect to $\mathcal G$.\par
In the case of $\mathcal G={\rm Diff(\mathbb R^3)}$, we obtain the following 
well-known fibration using corresponding symmetry group $G=SL(4,\mathbb R)$.
\setlength{\unitlength}{0.7mm}
\begin{equation}
\begin{picture}(100,50)(-30,0)
\put(-18,0){\makebox(20,10)[r]{$RP^3\cong SL(4,\mathbb R)/H_4$}}
\put(5,30){\makebox(20,10){$SL(4,\mathbb R)/H\cong \mathbb F_{V}(1,3)$}}
\put(28,0){\makebox(20,10)[l]{$SL(4,\mathbb R)/\overline H_1\cong RP^3$}}
\put(-34,17){\makebox(5,10){$RP^2\cong H_4/H$}}
\put(37,16){\makebox(35,10){$\overline H_1/H\cong RP^2$}}
\put(1,30){\vector(-1,-1){20}}
\put(20,30){\vector(1,-1){20}}
\end{picture}
\end{equation}
This fibration equals the fibration (14).\par
In the case of $\mathcal G$=ScaleDiff($\mathbb R^3$), we calculate the 
corresponding flat model space. 
From the action (15) of $SL(4, \mathbb R)$ on $\mathbb R^3$, we have the following 
symmetry group from restriction of variables associated with the scale transformation.\\
\begin{equation*}
G=\left\{
\begin{pmatrix}
*&0&0&* \\
0&*&0&* \\
*&*&*&* \\
0&0&0&*
\end{pmatrix}
\in SL(4,\mathbb R) \right\}. 
\end{equation*}
Then, we have the following.
\begin{proposition}
We obtain the following double fibration as the flat model space 
associated with $\mathcal G$={\rm ScaleDiff}$(\mathbb R^3)$:
\begin{equation}
\setlength{\unitlength}{0.7mm}
\begin{picture}(100,50)(-30,0)
\put(3,0){\makebox(-20,10)[r]{$G/(G\cap H_4)\cong \mathbb R^3$}}
\put(16,30){\makebox(-10,10){$G/(G\cap H)\cong \mathbb R^3$}}
\put(39,0){\makebox(20,10)[l]{$G/(G\cap \overline H_1)\cong \mathbb R$}}
\put(-20,17){\makebox(-50,10){$(G\cap H_4)/(G\cap H)\cong \left\{0\right\}$}}
\put(50,16){\makebox(30,10){$(G\cap \overline H_1)/(G\cap H)\cong \mathbb R^2$}}
\put(1,30){\vector(-1,-1){20}}
\put(20,30){\vector(1,-1){20}}
\end{picture}
\end{equation}
\begin{proof}
We prove the correspondence $G/(G\cap H_4)\cong \mathbb R^3$. 
We consider the following injective group homomorphism $\Phi:\mathbb R^3\to G/(G\cap H_4)$
 defined by:
$$
\Phi(a):=\left[
\begin{pmatrix}
1&0&0&a_1 \\
0&1&0&a_2 \\
0&0&1&a_3 \\
0&0&0&1
\end{pmatrix}
\right],$$
where, $a=(a_1,\ a_2,\ a_3) \in \mathbb R^3$. 
It is clear that $\Phi$ is bijective. Similarly, we have the other correspondences. Thus, we complete the proof.
\end{proof}
\end{proposition}

This fibration is degenerate. Since $\mathcal G={\rm ScaleDiff}(\mathbb R^3)$ is very strongly restricted 
from Diff($\mathbb R^3$), this degeneration arises. Hence, we 
consider the following problem.
\begin{problem}
Find a symmetry group for proper subgroup $\mathcal G$ of $\rm Diff(\mathbb R^3)$, from
which has double fibration of compact-type as a flat model space.
\end{problem}
 To consider this problem, we characterize groups $H_4,\ \overline H_1,\ H$. 
To this purpose, we prepare some terminology and notation. We use the 
Iwasawa decomposition:  
\begin{center}\label{iwasawa}
$SL(4, \mathbb R)=KAN,$ 
\end{center}
where, $K=SO(4)$ and $A=\left\{diag(a_1,a_2,a_3,a_4) \in SL(4,\mathbb R) \right\}$ and 
$N$ is a group of upper triangle matrices whose diagonal components are all 1. By 
using the Iwasawa decomposition, we have decompositions of the subgroups 
$H_4,\ \overline H_1,\ H$: 
\begin{center}
$H_4=K_4M_4A\overline N, \quad
\overline H_1=K_1M_1A\overline N, \quad
H=(K_1\cap K_4)M_{1, 4}A\overline N,$
\end{center}
where, $K_i\cong SO(3)$ are isotropy subgroups of $e_i \in S^3$ and 
$\overline N=\tilde\theta(N)$
and $M_{1}$, $M_{4}$ are the following subgroups of $M$: 
\begin{center}
  $M_{1}:=\langle1, m_1\rangle \cong \mathbb Z_2,\quad 
  M_{4}:=\langle1, m_4\rangle\cong \mathbb Z_2,\quad
  M_{1,4}:=M_{1}\cap M_{4},$
\end{center}
where, $m_1=diag(-1,-1,1,1)$ and $m_4=diag(1,1,-1,-1)$.

By using these facts, we find the fibration of compact-type. 
We consider the following subgroup $G \subset SL(4,\mathbb R)$:
\begin{align}\label{compact}
G&=\left\{g\in SL(4,\mathbb R) \ |\ g[e_3]=[e_3],\ ^{t}g^{-1}[e_3]=[e_3] \right\}, \\
 &=\left\{
\begin{pmatrix}
*&*&0&* \\
*&*&0&* \\
0&0&*&0 \\
*&*&0&* 
\end{pmatrix}
\in SL(4, \mathbb R) \right\}.
\end{align}
Note that $G$ is a subgroup invariant under Cartan involution $\tilde\theta$. We show 
that a double fibration defined by the group $G$ is a fibration of compact-type. 
To prove this assertion, we consider the characterization of groups $G,\ G\cap H_4,\
G\cap \overline H_1,\ G\cap H$. In fact, these groups are decomposited as follows.
\begin{lemma}
The above Lie groups $G,\ G\cap H_4,\ G\cap \overline H_1,\ G\cap H$ have the following 
decompositions:
\begin{align}
G&=(K\cap G)A(\overline N\cap G) \\
G\cap H_4&=(K_4\cap G)M_4A(\overline N\cap G) \\
G\cap \overline H_1&=(K_1\cap G)M_1A(\overline N\cap G) \\
G\cap H&=(K_1\cap K_4\cap G)M_{1,4}A(\overline N\cap G)
\end{align}
\begin{proof}
We prove (21). Since $K\cap G$, $A$, 
$\overline N\cap G$ are subgroups of $G$, we have\\ 
$G \supset (K\cap G)A(\overline N\cap G)$. 
Hence we prove $G \subset (K\cap G)A(\overline N\cap G)$. For $g \in G$, we write 
$g=kan$ ($k\in SO(4),\ a\in A,\ n\in \overline N$). By definition (19) of $G$, we assume as
follows:
\begin{center}
$ge_3=\alpha e_3, \quad
^{t}g^{-1}e_3=\gamma e_3,$
\end{center} 
where, $\alpha, \gamma \in \mathbb R^{*}$. 
We write $n \in \overline N$ and $a\in A$ explicitly as follows:
\begin{equation}
n=
\begin{pmatrix}
1&0&0&0 \\
n_1&1&0&0 \\
n_2&n_3&1&0 \\
n_4&n_5&n_6&1
\end{pmatrix},
\hspace{2cm}
a=
\begin{pmatrix}
a_1&0&0&0 \\
0&a_2&0&0 \\
0&0&a_3&0 \\
0&0&0&a_4
\end{pmatrix}.
\end{equation}
Then, 
\begin{align}
\alpha e_3&=kane_3=a_3ke_3+a_4n_6ke_4, \nonumber \\
 ^{t}ke_3&=\alpha^{-1}a_3e_3+\alpha^{-1}a_4n_6e_4.
\end{align}
On the other hand, if we write 
\begin{center}
$^{t}n^{-1}e_3=l_{2}e_1+l_3e_2+e_3, $
\end{center}
then, 
\begin{center}
$e_3=l_2e_1+l_3(n_1e_1+e_2)+n_2e_1+n_3e_2+e_3$. 
\end{center}
Hence, we have the equalities $l_3=-n_3,\ l_2=n_1n_3-n_2$ by using a linear independence of $e_1,\ e_2,\ e_3,\ e_4$. 
By using these relations, 
\begin{align}
^{t}g^{-1}e_3&=ka_1^{-1}(n_1n_3-n_2)e_1-ka_2^{-1}n_3e_2+ka_3^{-1}e_3, \nonumber \\
^{t}ke_3&=\gamma^{-1}a_1^{-1}(n_1n_3-n_2)e_1-\gamma^{-1}a_2^{-1}n_3e_2+\gamma^{-1}a_3^{-1}e_3
\end{align}
In the equality between right sides of (26) and (27), we have $n_2=n_3=n_6=0$ 
. Hence, we have $n\in G$. 
Since $g,\ a,\ n \in G$, we have $k\in G$. Thus we complete the proof of (21).  
By similar method, we can prove the $(22),\ (23),\ (24)$.
\end{proof}
\end{lemma}
 By using these decomposition formula, we have the following fibration of compact-type:
\begin{theorem}
A double fibration constructed by the group $(\ref{compact})$ is the following \\
fibration of compact-type.
\begin{equation}
\setlength{\unitlength}{0.7mm}
\begin{picture}(100,50)(-30,0)
\put(3,0){\makebox(-20,10)[r]{$G/(G\cap H_4)\cong  RP^2$}}
\put(16,30){\makebox(-10,10){$G/(G\cap H)\cong \mathbb F_{V^3}(1,2)$}}
\put(39,0){\makebox(25,10)[l]{$G/(G\cap \overline H_1)\cong RP^2$}}
\put(-15,17){\makebox(-50,10){$(G\cap H_4)/(G\cap H)\cong S^1$}}
\put(50,16){\makebox(30,10){$(G\cap \overline H_1)/(G\cap H)\cong S^1$}}
\put(1,30){\vector(-1,-1){20}}
\put(20,30){\vector(1,-1){20}}
\end{picture}
\end{equation}
\begin{proof}
By a direct computation, we have 
\begin{align*}
G\cap K&=SO(3)\times \mathbb Z_2, \quad
G\cap K_1=SO(2)\times \mathbb Z_2, \\
G\cap K_4&=SO(2)\times \mathbb Z_2, \quad
G\cap K_1\cap K_4=\mathbb Z_2. 
\end{align*}
Thus, we obtain the statement by the following correspondence:
\begin{align*}
G/(G\cap H_4)&\cong RP^2, \quad 
G/(G\cap \overline H_1)\cong RP^2, \quad
G/(G\cap H)\cong \mathbb F_{V^3}(1,2), \\
(G\cap H_4)/(G\cap H)&
\cong S^1, \quad 
(G\cap \overline H_1)/(G\cap H)\cong S^1.
\end{align*}
\end{proof}
\end{theorem}
We note that the coordinate transformation group $\mathcal G$ corresponding to 
this symmetry group $(\ref{compact})$ is constructed by the transformations of the form:
$$X_1=X_1(x_1,x_2),\ 
X_2=X_2(x_1,x_2),\ 
Y=\frac{y}{A(x_1,x_2)}.$$
\section{The dual equations}
\hspace{4mm}
In this section, we compute explicitly the dual equations of the second order PDE (1) 
([17]).\par
First, we assume that solutions of (1) are written by three parameters 
$X_1, X_2, Y$ as follows:
\begin{equation}
y=h(x_1,x_2,X_1,X_2,Y)
\end{equation}
Then, we can choose a local coordinate ($X_1$, $X_2$, $Y$) on a solution space of (1). 
A family of solutions corresponding to the above solutions is given by moving parameters 
$X_{1},\ X_{2},\ Y$ in the solution space. 
Let $Y=Y(X_1, X_2)$ be the surface on the solution space. Then, $Y$ can be written as 
$Y(X_1,X_2)=g(X_1,X_2,x_1,x_2,y)$.
We calculate $Y_{X_1},Y_{X_2}$ by using this representation:\\
By taking a derivation of $y=h(x_1,x_2,X_1,X_2,Y(X_1,X_2))$, we have 
\begin{align*}
Y_{X_i}=-\frac{h_{X_i}}{h_Y}.
\end{align*}
From this fact, the dual equation of (1):
\begin{equation}
\frac{\partial ^{2}  Y}{\partial X_i \partial X_j}=F_{ij}(X_{1},X_{2},Y,Z_{1},Z_{2})
\end{equation}
are written:
\begin{align*}
\frac{\partial ^{2}  Y}{\partial X_i \partial X_j}&=-\frac{1}{h_Y^2}\bigl\{(h_{X_iX_j}+h_{X_iY}Y_{X_j})h_Y-h_{X_i}(h_{YX_j}
+h_{YY}Y_{X_j})\bigr\} \\
&=\frac{h_{X_i}h_{YX_j}-h_Yh_{X_iX_j}+Z_j(h_{X_i}h_{YY}-h_Yh_{X_iY})}{h_Y^2}
\end{align*}
In particular, we calculate the  dual equation of (12). 
Solutions of (12) are written as $y=X_1x_1+X_2x_2+Y$, $z_i\equiv y_{x_i}=X_i$. Hence, $Y_{X_i}=-x_i$ and we obtain the following equations by substituting this term into above dual equations:
\begin{equation}
\frac{\partial ^{2}  Y}{\partial X_i \partial X_j}=0
\end{equation}
This equation is the dual equation of the flat equations (12). 
Namely, the dual equation of the flat equation is also the flat equation. 
This fact is supported by the double fibration (13). 
In the double fibration (13), the solution space of the original flat equation (12)
is a right base space $\mathbb R^3$, and the solution space of the dual equation (31) is a 
left base space $\mathbb R^3$.

\end{document}